\documentclass[graybox]{svmult}
\usepackage{graphicx}        % standard LaTeX graphics tool
\usepackage{multicol}        % used for the two-column index
\usepackage[bottom]{footmisc}% places footnotes at page bottom

\usepackage[T1]{fontenc}
\usepackage{url,amsmath,amsfonts}
\usepackage[mediumspace,mediumqspace,squaren]{SIunits}

\newcommand{\Prob}{{\mathbb P}}
\newcommand{\R}{{\mathbb R}}
\renewcommand{\E}{{\mathbb E}}
\newcommand{\F}{{\mathbb F}}

\newcommand{\Z}{{\mathbb Z}}
\newcommand{\N}{{\mathbb N}}

\newcommand{\C}{{\mathbb C}}

\newcommand{\cK}{{\mathcal K}}

\renewcommand{\L}{{\mathbb L}}
\newcommand{\ind}{{\mathbb I}}

\begin{document}

\title*{3D image based stochastic micro-structure modelling of foams
  for simulating elasticity}
\titlerunning{3D image based micro-structure modelling of foams for simulating elasticity}
\author{Anne Jung and Claudia Redenbach and Katja Schladitz and Sarah Staub}

\institute{Anne Jung \at Universität des Saarlandes, Technical Mechanics, Saarbrücken, Germany, \email{anne.jung@mx.uni-saarland.de}
\and Claudia Redenbach \at Technische Universit\"at Kaiserslautern, Mathematics Department, Kaiserslautern, Germany,  \email{redenbach@mathematik.uni-kl.de}
\and Katja Schladitz \at Fraunhofer-Institut f\"ur Techno- und Wirtschaftsmathematik, Kaiserslautern, Germany, \email{katja.schladitz@itwm.fraunhofer.de}
\and Sarah Staub \at Fraunhofer-Institut f\"ur Techno- und Wirtschaftsmathematik, Kaiserslautern, Germany,  \email{sarah.staub@itwm.fraunhofer.de}
}

\maketitle

\abstract{Image acquisition techniques such as micro-computed
  tomography are nowadays widely available. Quantitative analysis of
  the resulting 3D image data enables geometric characterization of
  the micro-structure of materials. Stochastic geometry models can be
  fit to the observed micro-structures. By alteration of the model
  parameters, virtual micro-structures with modified geometry can be
  generated. Numerical simulation of elastic properties in
  realizations of these models yields deeper insight on the influence
  of particular micro-structural features. Ultimately, this allows for
  an optimization of the micro-structure geometry for particular
  applications.  Here, we present this workflow at the example of open
  cell foams. Applicability is demonstrated using an aluminum alloy
  foam sample. The structure observed in a micro-computed tomography
  image is modelled by the edge system of a random Laguerre
  tesselation generated by a system of closely packed spheres. Elastic
  moduli are computed in the binarized $\micro$CT image of the foam as
  well as in realizations of the model. They agree well with the
  results of a compression test on the real material.\footnote{This
    is a preprint of the following chapter: Anne Jung, Claudia
    Redenbach, Katja Schladitz, and Sarah Staub: 3D Image-Based
    Stochastic Micro-structure Modelling of Foams for Simulating
    Elasticity, published in Research in Mathematics of Materials
    Science. Association for Women in Mathematics Series, vol 31,
    edited by Malena I. Español, Marta Lewicka, Lucia Scardia, and
    Anja Schlömerkemper, 2022, Springer, Cham, reproduced with
    permission of Springer Nature Switzerland AG 2022.  The final
    authenticated version is available online at:
    {https://doi.org/10.1007/978-3-031-04496-0\_11}} 
}

\section{Introduction}
\label{sec:introduction}

3D image data, predominantly generated by computed tomography, have
been analyzed for more than 30 years now. In medical applications, the
emphasis has been on segmentation and visualization. Algorithmic
foundations had been settled in the 1990es, see \cite{lohm98} for a
comprehensive overview. 

Quantitative analysis of 3D images of micro-structures started in the 1990es, too, 
with porous rock samples \cite{biswal98,lindquist94} and trabecular bone \cite{hildebrand97}.
Prediction of macroscopic, in particular mechanical, properties based on the observed micro-structure has been a goal right from the beginning \cite{ruegsegger97}. 

The seminal paper \cite{lan:ohs:hil99} provided Ohser's algorithm for efficiently estimating the intrinsic volume densities based on 3D binary images,  
linking fundamental characteristics of random closed sets \cite{KMS2013,hadwiger57,schn:wei08} to those observable in image data. 

\begin{figure}[b]
\sidecaption
\includegraphics[width=.32\textwidth]{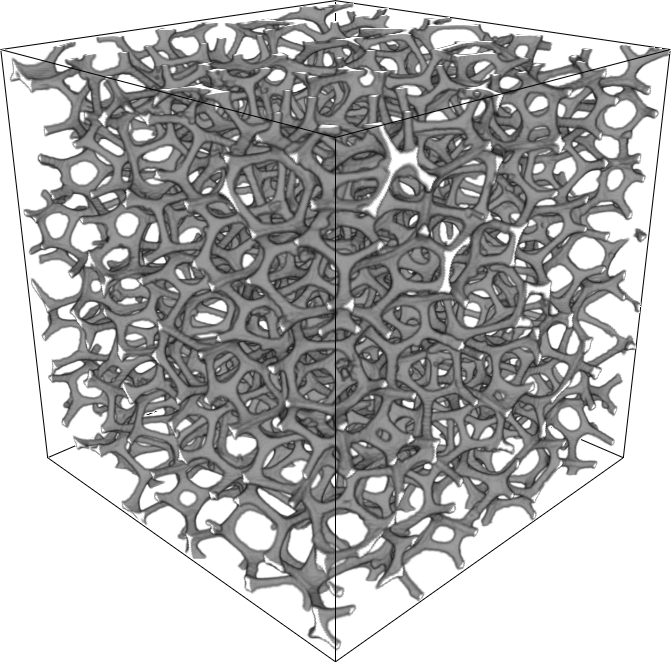}\
\includegraphics[width=.32\textwidth]{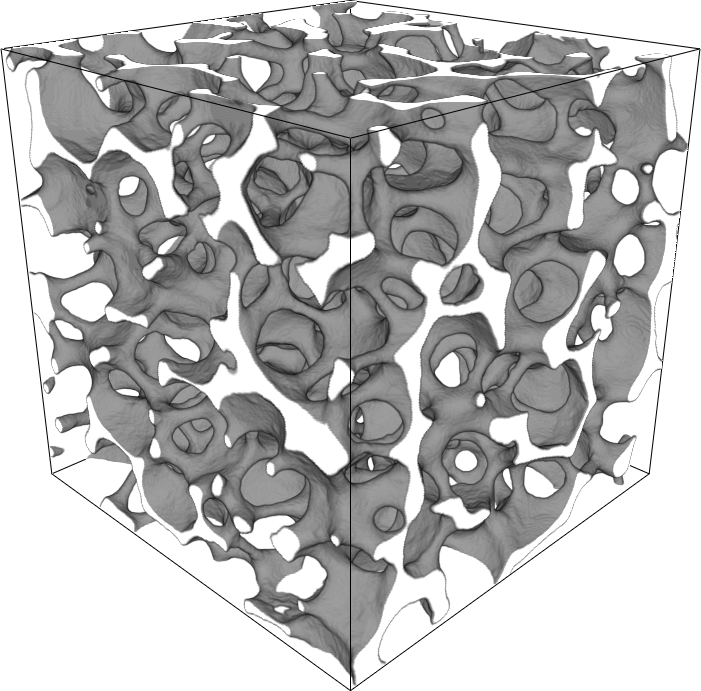}\
\includegraphics[width=.32\textwidth]{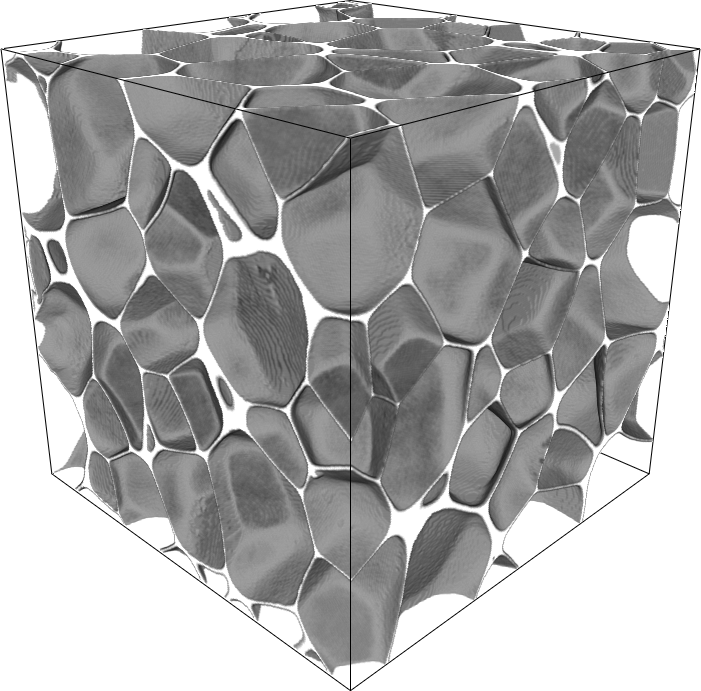}
\caption{Visualizations of $500^3$ voxel sub-volumes from reconstructed CT images of rigid foams. From left to right: open polymer, partially closed ceramic, and closed polymer foams. All CT scans taken at ITWM with voxel sizes $49\,\micro$m for the open polymer foam, $34\,\micro$m for the ceramic foam, and $3\,\micro$m for the closed polymer foam. Samples provided by Vesuvius (ceramic) and Evonik (closed polymer foam).}
\label{fig:severalFoams}
\end{figure}

Our contributions to predicting macroscopic materials properties by numerical simulation in realizations of 3D stochastic 
geometry models started with acoustic adsorption \cite{schladitz06:_desig} and had been concerned with rigid foams early on
\cite{helfen03:_inves_foamin_proces_metal_synch_radiat_imagin,rack08:_analy}. After establishing the mathematical 
\cite{lautensack06:_random_laguer} and algorithmic \cite{lautensack08:_fittin_laguer} bases for modelling cellular materials 
structures, elastic properties were studied \cite{redenbach2012}. Here, we will concentrate on characterization and modelling of rigid foams, too. Some typical examples of these materials are shown in Figure~\ref{fig:severalFoams}.

Stochastic modelling allows for detailed studies of structure property relations
\cite{foehst:redenbach21,liebscher2012,vecchio16}. Geometric features closely tied in practice can be modified 
individually. The use of stochastic instead of deterministic geometry models naturally captures microscopic 
heterogeneity as well as macroscopic homogeneity assuming invariances of the underlying distribution law. 
Moreover, the size of the representative volume element (RVE) can be determined 
statistically \cite{jeulin21,kanit03}. 
In this work, we investigate RVE of constant size corresponding to the size of the foam 
samples characterized experimentally. Alternatively, more, but smaller volume elements 
in the sense of stochastic homogenization as devised e.~g. in \cite{armstrong2019} could be used. 
In the case of rigid foams, the size of these volume elements is surely bounded by one foam cell. 
The theory of homogeneous random closed sets provides the corresponding mathematical concept by 
the typical cell of a tessellation \cite{schn:wei08}. Exploring this alternative approach and 
comparing the two are subject of future research.
Fitting stochastic geometry models to the observed real structure is tedious as analytic relations are 
available only for models not  suitable for real materials. 
Starting at the real micro-structure ensures however that observed trends and 
correlations have a practical meaning.

The concept of statistically similar representative volume elements (SSRVE) \cite{balzani2014} is similar to our approach in the sense of fitting the synthetic structure to the real one by minimizing the weighted sum of the squared differences of several geometric characteristics estimated for both, real and synthetic structure. It differs however fundamentally in aiming at deterministic, geometrically significantly simplified structures right from the beginning. Simple structures consisting e.g. of just a couple of spheres are used to generate the SSRVE that is subsequently periodically continued.

Here, we focus on numerical computation of effective elastic properties of an 
aluminum alloy foam by homogenization. The foam sample's micro-structure is observed in 3D image data obtained by micro-computed tomography ($\micro$CT). 
We use both, the segmented image data as well as a random Laguerre tessellation model fit to the observed foam structure based on the estimated geometric characteristics. 
We provide the basics on random closed sets, their characteristics, in particular the densities of the intrinsic volumes, and estimating them based on 3D image data in Section~\ref{sec:image-analysis}. 
The image processing workflow for dividing the pore space of open foams into individual image objects is described in the same section. Section~\ref{sec:laguerre-tessellations} is dedicated to random Laguerre tessellations and fitting them to the observed structure. 

Simulation of elastic properties in voxel representations of the real foam or of 
model realizations is described in Section~\ref{sec:Homogenization}.
A computational homogenization scheme \cite{Suquet1985} is used to
transition from the microscopic scale of the foam structure to the macroscopic scale of 
the effective properties. A comprehensive overview on this class of homogenization 
schemes is given in \cite{Geers2010}. A microscopic boundary value problem is formulated 
by imposing admissible boundary conditions and subsequently solved numerically,
by a Finite Element Method \cite{Feyel2000} or a fast Fourier transform (FFT) based 
solution of the Lippmann-Schwinger (LS) equations of elasticity \cite{Moulinec1998}. 
Efficiency in terms of computational effort and memory use is critical as the images consist of  several hundreds of voxels in each direction. The LS-FFT homogenization applied here is described in Section~\ref{sec:Homogenization}.

Finally, in Section~\ref{sec:application}, we apply the methods from 
Sections~\ref{sec:image-analysis}-\ref{sec:Homogenization} to the real aluminum foam 
sample. A random Laguerre tessellation is fit to the foam structure based on 
the characteristics estimated from the segmented CT image. Four synthetic foams 
are derived from its edge system by applying two cross-sectional shapes of the struts 
and relaxing the foam or not.
Elastic properties derived by computation in the segmented $\micro$CT image and the 
structural modulus measured experimentally are compared in \ref{subsec:mechprops}. The 
effective stiffnesses of the synthetic foams are numerically predicted in 
\ref{ssubsec:prediction}.

\section{3D image analysis for foams}
\label{sec:image-analysis}

We introduce the general concept of random closed sets, define basic characteristics for them, and describe how to estimate these characteristics based on 3D volume images as generated by computed tomography. 

\subsection{Random closed sets and their characteristics}
\label{subsec:random-closed-sets}

Materials micro-structures are often macroscopically homogeneous 
in some sense but locally heterogeneous.
A commonly used mathematical model for such structures are 
random closed sets, random
variables whose realizations are closed subsets of $\R^3$. 
Schneider and Weil \cite{schn:wei08} attribute the concept to Matheron \cite{matheron72racs} and Kendall \cite{kendall74racs}, with the first detailed description being \cite{mat75}.

Denote by ${\mathcal F}$ the system of closed subsets of $\R^3$ and 
by $\mathfrak F$ the hit-or-miss $\sigma$-algebra 
generated by the sets $\{F\in\mathcal F: F\cap A\neq\emptyset\}$ for
all compact $A\subset\R^3$. The pair $({\mathcal F},\mathfrak F)$ is 
then a measurable space and we can define random closed sets (RACS) as 
random variables with values in this space. More precisely, 
let $(\Omega, \mathcal A, \Prob)$ be a probability space. 
A random closed set $\Xi$ is a measurable mapping
\[\Xi: (\Omega, \mathcal A, \Prob)\mapsto({\mathcal F},\mathfrak F).\] 
We call $\Xi$ stationary (or macroscopically
homogeneous) and isotropic if its distribution is invariant with respect to translations and rotations, respectively. 

\paragraph{Point processes}

One class of RACS models are locally finite random closed sets: 
Denote by
${\mathcal C}$ the system of compact subsets of $\R^3$. A set $A\subseteq\R^3$ is 
locally finite if $\#(A \cap C) <\infty$ for any $C \in {\mathcal C}$. Here $\# B$ refers to the number of elements of the set B. A (simple)
point process $\Phi$ is a random variable with values in
the system of locally finite sets and can be read as a sequence of random points in $\R^3$ as well as a random counting measure. 
The measure $\Lambda$ defined by $\Lambda(B)=\E(\#(\Phi \cap B))$ is called
the intensity measure of $\Phi$. For stationary point
processes, it is a multiple of the Lebesgue measure 
$\Lambda(B) = \lambda V(B)$ with the intensity $\lambda >0$ giving the
mean number of points per unit volume.

\paragraph{The Poisson point process}

The Poisson point process representing complete spatial randomness plays a
particular role as it is analytically tractable and the basis for 
constructing other point process as well as a variety of random closed set 
models. A Poisson point process $\Phi$ in $\R^3$ with intensity measure 
$\Lambda$ is a point process with the number of points $\#(\Phi \cap B)$ in a Borel set $B$ being Poisson distributed with parameter $\Lambda(B)$ and 
the numbers of points $\#(\Phi
\cap B_1), \ldots, \#(\Phi \cap B_n)$ in pairwise disjoint 
Borel sets $B_1, \ldots, B_n $ being independent random variables.

As a consequence, a stationary Poisson point process with intensity $\lambda$ in
a compact observation window $W\subset\R^3$ can easily be simulated by first drawing the number of points Poisson distributed with parameter $\lambda V(W)$. Then,
the points are placed
independently identically uniformly distributed in $W$. There is no interaction of the points. This is why the Poisson
point process is used as reference for complete spatial randomness.

The random points can additionally be decorated by marks representing for 
instance a species, age, or size. Moreover, the definitions of point process 
and Poisson point process can directly be transferred any locally compact space 
with countable base.

\subsubsection{Intrinsic volumes}
\label{subsubsec:intrinsic-volumes}

The intrinsic volumes are basic characteristics of RACS. 
Let $\cK$ denote the system of compact and convex 
sets (convex bodies) in $\R^3$.  Dilating $K\in \cK$ with a ball
of radius $r$ yields the so-called parallel set 
$K \oplus b(0,r) := \{x+y\,:\, x \in K, y \in b(0,r)\}$.  
Steiner's formula expresses the volume of this set as a polynomial in $r$:
\[V(K \oplus b(0,r)) = \sum_{k=0}^3 \kappa_k V_{3-k}(K) r^k, \quad
r\ge 0,\] 
where $\kappa_k$ is the volume of the $k$-dimensional unit
ball and the coefficients $V_0,\ldots, V_3$, are the four intrinsic
volumes in $\R^3$. These are -- up to constant factors
-- the volume $V=V_3$, the surface area $S=2V_2$, the integral of mean
curvature $M=\pi V_1$, and the Euler number $\chi=V_0$, see
e.~g. \cite[p. 210]{Schn93}. For $K\in \cK$, the integral of mean curvature is up to a constant the
mean width $M= 2\pi{\bar{b}}$ -- the distance of two parallel planes enclosing
$K$, averaged w.r.t. rotation. 

The intrinsic volumes form a basis of geometric characteristics in the sense of 
Hadwiger's theorem 
\cite[Theorem 14.4.6]{schn:wei08}: Every rigid motion invariant, additive, and 
continuous functional on $\cK$ is a linear combination of the intrinsic volumes. 
Additivity yields a straightforward extension of the intrinsic volumes to finite unions of convex bodies, so-called polyconvex sets.

\subsubsection{Densities of the intrinsic volumes}
\label{subsubsec:densities-intrinsic-volumes}

The intrinsic volumes can be used to characterize single cells of a foam. For characterizing the complete solid component, we need another generalization. %to the case of observing a subsample of a stationary RACS in a compact convex observation window. 
It applies for RACS $\Xi$ whose intersections with convex, compact observation windows $W \in \cK$ are almost surely polyconvex. 
In this case, the intrinsic volumes of $\Xi \cap W$ are well-defined, too. 
Denote by $\E V_k(\Xi\cap W)$ the expected $k$-volumes of this intersection with respect 
to the distribution of $\Xi$. 

The {\em densities of the intrinsic volumes} are the limits of these expectations for growing observation window: 
\begin{equation} 
\label{eq:intrinsic-volume-densities}
V_{V,k}(\Xi)=\lim\limits_{r\to \infty}
\frac{\E V_k(\Xi\cap rW)}{V(rW)},\qquad k=0,\ldots , 3.
\end{equation}
The limit exists if $\E 2^{N(\Xi \cap [0,1]^3)}< \infty$ \cite[Theorem 9.2.1]{schn:wei08}, where $N(X)$
is the smallest $m\in\N$ such that $X= K_1 \cup \ldots \cup K_m$ with 
$K_1, \ldots, K_m \in \cK$. 

Of particular interest for foams are the (solid) volume fraction $V_V=V_{V,3}$, with $p= 1-V_V$ being the porosity, and the specific surface area $S_V= 2 V_{V,2}$. Additionally, the density of the Euler characteristic $\chi_V= V_{V,0}$ is related to the mean number of nodes of the foam skeleton per unit volume \cite{open-foam-features-2011}.

\subsection{Image analysis}
\label{subsec:image-analysis}
We quickly define the basic concepts image and adjacency system, discuss the segmentation tasks 
posed by our problem, and finally describe how to estimate the characteristics needed for model fitting 
later on.

\subsubsection{Images}
\label{subsubsec:images}

A typical way of observing a realization $X$ of a
suitable RACS $\Xi$ is via discretization in a 3D image. For simplicity, let
$\L=s\Z^3$ be a three-dimensional cubic lattice with
lattice spacing $s>0$. Denote by $W\subset\R^3$ a cuboidal observation
window. By an image we understand a function
\[f:\L\cap W\longrightarrow V,\] where $V$ is the set of real or
complex numbers $\R$, $\C$, or $V=\{0,\dots,2^n-1\}$ with $n=1,8,16$
or $32$. 

Black-and-white images $f$ with $V=\{0,1\}$ are often called binary image. 
In this case, the image foreground 
can be identified with the intersection
$X\cap\L\cap W$ of a set $X\subset\R^3$ observed at the
lattice points in the observation window $W$. The function $f$ is then the indicator
function of $X$, restricted to the observable points
$f=\ind_X|_{\L\cap W}$. That is, $f(x)=1$ if $x\in X\cap\L\cap W$
and $f(x)=0$ otherwise. 
% Binary images can be derived from real-valued ones, and we will give
% a brief introduction into this topic in Sec.~\ref{sec:segmentation}.
Usually, in image processing, both
$x\in\L\cap W$ and the pair $(x,f(x))$ are called pixel (or voxel). 

Mathematically, we interpret an image as a set of lattice points 
equipped with grayvalues, not as a set of small cubes. We nevertheless 
address the lattice spacing $s$ as pixel or voxel size, too.

By construction, the pixels of a 3D image are discrete. To represent concepts such as the connected components of a RACS in an image, a notion of connectivity of pixels must be introduced. To define the discrete
connectivity in 3D unambiguously, we follow \cite{oh:na:sc02} in using
adjacency systems consisting of $j$-dimensional faces with
$j=0,\dots,3$. The \emph{discretization} $X\sqcap\F$ of a compact subset
$X\subset\R^3$ w.~r.~t. a given adjacency system $\F$ is defined as the
union of the elements $F$ of $\F$ whose vertices ${\mathcal F}^0(F)$ are 
contained in $X$
\[X\sqcap\F=\bigcup\{F\in\F\,:\,{\mathcal F}^0(F)\subseteq X\}.\]
It can be interpreted as an approximation of $X$ by a polyhedral set
built using the bricks provided by the adjacency system. 
For background on adjacency systems and consistent connectivity for
foreground and background of an image, see 
\cite{nag:ohs:pis99,oh:na:sc02,ohs:nag:sch03,ohser-schladitz09book,schladitz06a}. 

\subsubsection{Segmentation}
\label{ssec:segmentation}

Segmentation subsumes a wide variety of methods assigning pixel values to classes, usually coded by integer values. Here, we have to solve two segmentation tasks: First, we have to binarize the image $f$ to get hold of the solid foam structure. Second, we exploit moments of the distributions of cell characteristics in the model fitting step. Hence, we have to identify foam  cells or pores as image objects. This second segmentation task is much more demanding, as the cells of an open foam are connected. As a consequence, the background or pore space has to be divided to "reconstruct" the cell structure.

\paragraph{Binarization}
\label{par:binarization}

Very roughly, computed tomography generates 3D images with 
voxel grayvalues related to the mean X-ray absorption  of the small sub-volume they represent. Solid material structure thus appears brighter than air. For binarizing the image, we look for a transformation 
\[T(f):\L\cap W\longrightarrow \{0,1\}\] 
such that $T(f)\equiv1$ corresponds to the solid foam material 
$X$ observed in $\L \cap W$ as described above.  
The simplest way to achieve a binarization, is just to use a global grayvalue threshold $\theta\in\R$
\[T(f)(x)=\begin{cases}
    1 & f(x) \geq \theta \\ 
    0 & \mathrm{else} \\
          \end{cases}\quad .
\]
This simple binarization method can of course be applied only, 
if the image is free of global grayvalue fluctuations. 
Prior to binarization, we apply a median filter with a $3\times3\times3$ filter mask to reduce noise and smooth the foam's surface slightly. 
Subsequently, the global grayvalue threshold $\theta$ can be chosen by 
Otsu's method \cite{otsu79}.

\paragraph{Cell reconstruction}
\label{par:cell-rec}

Now we separate the multiply connected pore system of the foam 
into individual image objects 
by a combination of two strong morphological tools, 
the Euclidean distance transform and the watershed transform.  

Let $X\subseteq\R^3$ be the set
under consideration, here the pore space of the foam. The Euclidean distance transform maps each point
in $\R^3$ to its shortest distance to the complementary set
$\R^3\setminus X$,
\[{\rm EDT}:\R^3\mapsto [0,\infty):x\mapsto\min\{||x-y||:
y\in\R^3\setminus X\}.\] 
This results in ${\rm EDT}(y)=0$ for all $y\in\R^3\setminus X$ and local maxima 
in the centers of spherical regions in $X$. Inverting the EDT image ${\rm
  EDT}(X\cap\L\cap W)$ turns these local maxima into minima:
$f(x)=\max\{{\rm EDT}(X\cap\L\cap W)\}-{\rm EDT}(x)$. The exact Euclidean distance 
transforms can be calculated very efficiently, i.e. in linear time, exploiting the Voronoi paradigm \cite{mau:rag03}. 

Now the watershed transform assigns a connected region to each local
minimum. This transform can be interpreted as the flooding of the
topographic surface $\{(x,f(x)): x\in\L\cap W \}$: water rises
uniformly with growing grayvalue $f$ from all local
minima. Watersheds are formed by all pixels where water basins
filled from different sources meet. Finally,
the image is segmented into regions and the system of watersheds
dividing them \cite{vin:soi91}. 

\begin{figure}[b]
\sidecaption
\includegraphics[width=.24\textwidth]{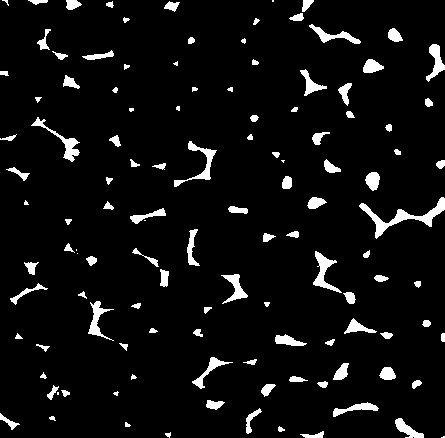}
\hfill
\includegraphics[width=.24\textwidth]{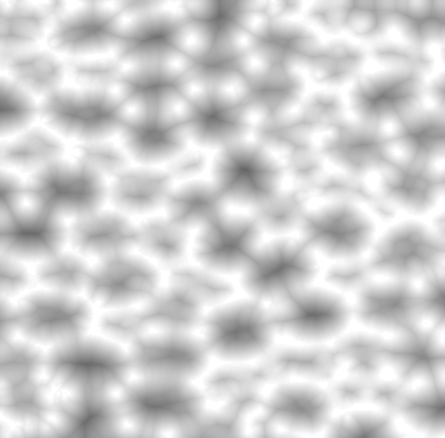}
\hfill
\includegraphics[width=.24\textwidth]{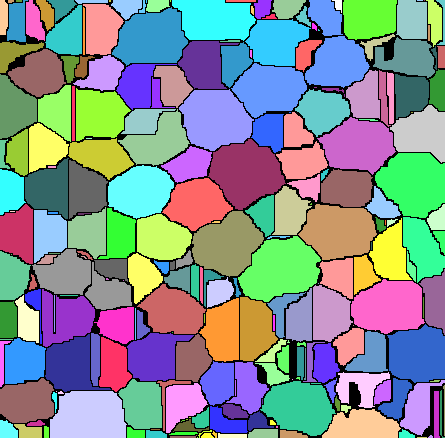}
\hfill
\includegraphics[width=.24\textwidth]{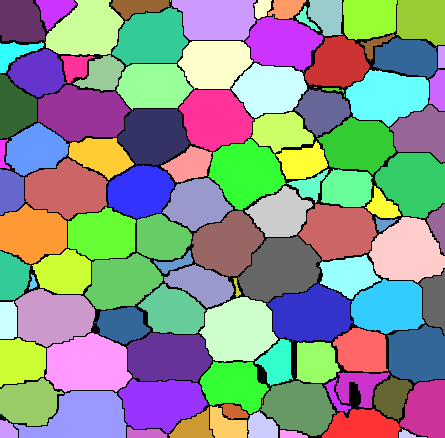}
\caption{Morphological cell reconstruction, illustrated using 
2D slices of the 3D image. Left: Binarized - solid foam structure appearing white, pore space black. Center, left: Inverted Euclidean distance map on the pore space, small values dark, high white. Center, right: Pore system generated by the watershed transform. Right: Pore system generated by the pre-flooded watershed transform.}
\label{fig:cell-rec}      
\end{figure}

In practice, this morphological separation strategy suffers from superfluous local minima in the inverted EDT (iEDT) image due to discretization, imperfect binarization, and cell shapes deviating from spherical. The watershed transform assigns each such minimum a region, inevitably resulting in a 
strong over-segmentation. A remedy for this is to prevent small regions from 
arising at all by altering the flooding accordingly (pre-flooded watershed, \cite{tek.dempster.ea:blood}). 
The foam structure under consideration here is rather regular. Predescribing a minimal 
cell volume is thus easily applicable. See Figure~\ref{fig:cell-rec} for an illustration using 2D sections 
through the 3D images.

\begin{svgraybox}
	It is remarkable that the heavyweights in this segmentation procedure -- watershed transform and EDT -- as devised by \cite{vin:soi91} and \cite{mau:rag03} work "as is" in 
	arbitrary dimensions. Often, image processing methods are mathematically easily formulated 
	for arbitrary dimensions and applied to 2D images. Practical application in higher dimensions is 
	nevertheless often impossible due to computational complexity or ambiguities, e.g. arising from the existence 
	of two subdimensions or connectivity issues.
\end{svgraybox}

% regions smaller than a pre-defined volume are united with neighbour regions
%\cite{tek.dempster.ea:blood} or smoothing the EDT image using the
%morphological h-minima transform \cite[Chapter 6]{soi99} are
%well-tested methods for this purpose. The latter even allows to cope with
%large size differences when the parameter h is adapted to total grey
%value \cite{godehardt06:_geomet_charac_light_weigh_compos},
%\cite[Section 4.2.6]{ohser-schladitz09book}.

\subsubsection{Estimating the intrinsic volumes based on image data}
\label{subsubsec:est-intrinsic-vol}

The intrinsic volumes and their densities can be estimated 
based on 3D image data 
very efficiently by Ohser’s algorithm \cite{lan:ohs:hil99}. 
This algorithm codes the 
$2\times2\times2$ voxel configurations in a binary image using a 
convolution with a cubic 
filter mask assigning weights $2^k,\ k=0,\ldots,7$ to the vertices of the cube. 
The gray value histogram of the 
resulting 8bit gray value image contains all information needed to 
derive the intrinsic volumes by 
multiplication with a suitable weights vector. 
See \cite{ohser-schladitz09book} for details.

\section{Random Laguerre tessellations and fitting them}
\label{sec:laguerre-tessellations}

Laguerre tessellations are a generalization of the well-known Voronoi tessellation model.  When generated by random sphere packings, they reflect the topology of rigid foams very well and enable particularly good control over the volume distribution of the resulting cells. We recall the concept and a fitting strategy. 
%Particle processes are a very versatile and general class of random set models  consisting of a set of randomly placed particles of random
%shape and orientation. Random tessellations are special particle processes whose
%particles are bounded cells dividing the $\R^3$. 
%They are thus well suited for modelling the micro-structure of cellular and 
%polycrystalline materials \cite{}. 

\subsection{Laguerre tessellations generated by random sphere packings}
\label{subsec:laguerre-tessellations}

Let $T$ be a set of bounded convex three-dimensional
subsets of $\R^3$, the cells of $T$. The system $T$ is called a tessellation of $\R^3$ if $T$ is space-filling, i.e. $\bigcup_{C\in T}C=\R^3$, and if the interiors of different cells do not intersect. Convexity of the cells and $T$ being space-filling force the cells to be threedimensional polytopes (3-polytopes) \cite[Lemma 10.1.1]{schn:wei08}. 

Write
\[F(x)=\bigcap_{C\in T,x\in C} C,\quad x\in\R^3\] 
for the intersection of all cells of $T$ containing the point
$x$. Then $F(x)$ is a non-empty finite intersection of 3-polytopes
and hence a $k$-polytope, $k\in\{0,\dots,3\}$.  The set
\[\Delta^k(T):=\{F(x):\dim F(x)=k,\quad x\in\R^3\},\quad k=0,\dots,3,\]
is the set of $k$-faces of the  tessellation $T$. We denote by
$\mathcal{F}^k(C)$ the set of all $k$-faces, $k=0,\dots,3$,
of a 3-polytope $C$. For a set $\F$ of convex polytopes write 
${\mathcal F}^k(\F)=\bigcup \{ {\mathcal F}^k(F): F\in\F \}$.

The tessellation $T$ is called face-to-face if the faces of the 
tessellation coincide with the faces of the cells. If additionally
exactly four cells meet in a vertex and three cells meet in an edge, the tessellation is called normal. 
A random tessellation is face-to-face or normal if its realizations almost surely have these properties.

Due to Plateau's laws, soap froths are normal, too \cite{book:WeHu1999}. 
Normal tessellations are thus a natural choice for modelling foam structures. 
The best known and most used tessellation model is the Voronoi 
tessellation that can be defined in arbitrary dimension $n$ based on a locally finite set $\varphi\subset\R^n$ of
generators. The Voronoi tessellation of $\varphi$ consists of the cells
\[C\big(x):= \{z \in \R^n \, :\, || x-z||\le ||y-z|| , \mbox{ for all }
y\in\varphi\}, \ x \in \varphi.\] 
The size of a Voronoi cell thus only depends on the distance of its generator to the neighboring generators.

To gain more flexibility regarding cell shapes and in particular sizes, we 
use Laguerre tessellations. In this weighted generalization of the 
Voronoi tessellation, each generating point $x \in \R^n$ is assigned 
a positive weight $r>0$ that can be interpreted as radius of a sphere with center $x$. 
The Laguerre cell of $(x,r)\in\varphi$ is defined as
\[C\big((x,r),\varphi\big):= \{z \in \R^3 \, :\, 
  ||x-z||^2 -r^2\le||y-z||^2-s^2,\ \mbox{ for all } (y,s)\in\varphi\}.\] 

The Laguerre tessellation $L(\varphi)$ is the
set of the non-empty Laguerre cells of $\varphi$. 
The Voronoi tessellation corresponds to the special case of the Laguerre tessellation generated by a system of spheres of constant radius. 

Laguerre tessellations are the most general model satisfying our 
assumptions: Each normal tessellation of $\R^3$ with convex cells is a
Laguerre tessellation \cite{Aurenhammer87a,LauZuy08}. 

\begin{figure}[b]
\sidecaption
\includegraphics[width=.32\textwidth]{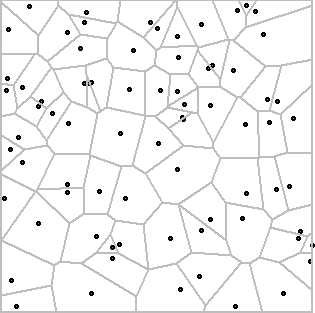}\
\includegraphics[width=.32\textwidth]{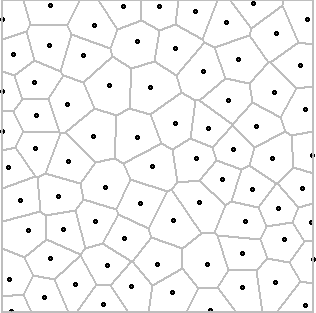}\
\includegraphics[width=.32\textwidth]{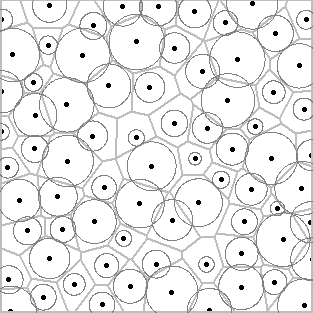}
\caption{Left and middle: Realizations of Voronoi tessellations generated by a Poisson point process and a regular point process. Right: Laguerre tessellation generated by the same point pattern as in the middle.}
\label{fig:VoronoiLaguerre}
\end{figure}

In Laguerre tessellations, empty cells as well as cells not containing
their generators can occur. Laguerre tessellations generated by 
non-overlapping spheres are however free of such irregularities, see Figure~\ref{fig:VoronoiLaguerre} for examples in $\R^2$. 

Random systems of closely packed spheres like those generated by the 
force-biased algorithm
\cite{bezrukov01,bezrukov02} are well suited to generate 
Laguerre tessellations with prescribed cell volume distribution \cite{Redenbach2009}. The 
cell shape is however rather restricted. In particular, the edge length distributions in Laguerre tessellations generated by dense packings of spheres are known to differ significantly from those observed in real foams, see \cite{vecchio16}. This problem can be alleviated by relaxing the structure using the Surface Evolver \cite{Brakke1992}.  

Voronoi or Laguerre tessellations generated by stationary and isotropic (marked) point processes are also isotropic. In contrast, real foam structures often show anisotropies as their cells are elongated in particular directions due to the foam generation process. A simple way of incorporating this anisotropy into the models is by appropriate scaling of the cell systems.

\subsection{Fitting a tessellation model}
\label{subsec:modelfitting}
Two strategies of fitting Laguerre tessellation models can be found in the literature. Laguerre approximation aims at finding a Laguerre tessellation that represents the observed cell system best in the sense that discrepancy between the observed cells and the cells of the approximation is minimized \cite{LiebscherApprox, Spettl}. Alternatively, a parametric tessellation model can be fit. That is, the observed cell system is approximated in a stochastic sense. The model is supposed to fit distributions of cell characteristics such as the volume, surface area or number of facets \cite{lautensack08:_fittin_laguer,Redenbach2009}. The former approach yields an exact representation of the observed structure in the class of Laguerre tessellations. Here, we concentrate on the latter approach, as it allows for generating an arbitrary number of model realizations of basically arbitrary size. 

The deviation of the model realization from the real observed foam 
structure is measured using the relative distance measure
\begin{equation}
\label{redenbach_eqn_dist}
\rho(\hat{m},m)=\sqrt{\sum_{i=1}^n \Big( \frac{m_i-\hat{m}_i}{\hat{m}_i} \Big)^2}, 
\end{equation}
where $\hat{m}=(\hat{m}_1,\ldots, \hat{m}_n)$ and
$m=(m_1,\ldots, m_n)$ are moments of the distributions of geometric characteristics
of the cells of the original foam and the model, respectively. 
For rigid foams whose structure is observed in 3D images, 
the means and standard deviations of
the volume $V=V_3$, the surface area $S=2 V_2$, the mean width $\bar{b}= \frac{1}{2} V_1$, and
the number of facets $F_C$ of the cells have proven to be particularly 
well suited \cite{lautensack08:_fittin_laguer,Redenbach2009}. The choice is based on mean value relations for normal random tessellations, see \cite[Section 9.4]{KMS2013}.

For minimizing the distance measure \eqref{redenbach_eqn_dist} on the parameter space of a tessellation model, one would ideally use analytic formulas relating the model parameters and the required moments of cell characteristics. Unfortunately, only Laguerre tessellations generated by a homogeneous Poisson process are analytically tractable \cite{LauZuy08}. For Laguerre tessellations generated by random sphere packings, model fitting therefore has to rely on Monte Carlo simulations of model realizations.

The cell volumes in cellular materials are usually assumed to be
lognormal or gamma distributed. Moreover, Laguerre
tessellations of dense packings of spheres with lognormal and gamma
distributed volumes are very regular and therefore well suited to model 
rigid foams. In \cite{Redenbach2009}, model realizations for
various packing fractions $\kappa$ and coefficients of variation $c$
of the volume distribution were generated. Subsequently, polynomials
in $c$ were fitted to the estimated geometric characteristics for each
value of $\kappa$. Using these results, the minimization of 
$\rho(\hat{m}, m)$ in \eqref{redenbach_eqn_dist} reduces to
minimizing a polynomial thus allowing for quick and easy model
fit.

\section{Numerical simulation of elastic properties}
\label{sec:Homogenization}
 The following subsections outline the simulation of the effective mechanical stiffness of inhomogeneous micro-structures such as foams.
 Therefore, the averaged stresses and strains are defined as volume averages over their microscopic counterparts.
 Furthermore, the fundamentals of the solution of the microscopic boundary value problem in terms of an LS-FFT are outlined briefly.
 
\subsection{Effective properties of micro-structured materials}
\label{ssec:Effective}

We apply a computational homogenization scheme to determine the effective stiffness of the foam numerically.
In this scheme, the microscopic geometry of the structure is captured either by a CT image of the foam or a realization of a stochastic geometry model.
Furthermore, the mechanical properties of the micro-constituents, here the  Young's modulus and Poisson's ratio, are required.

\begin{figure}[ht!]
\sidecaption
\includegraphics[width=0.95\textwidth]{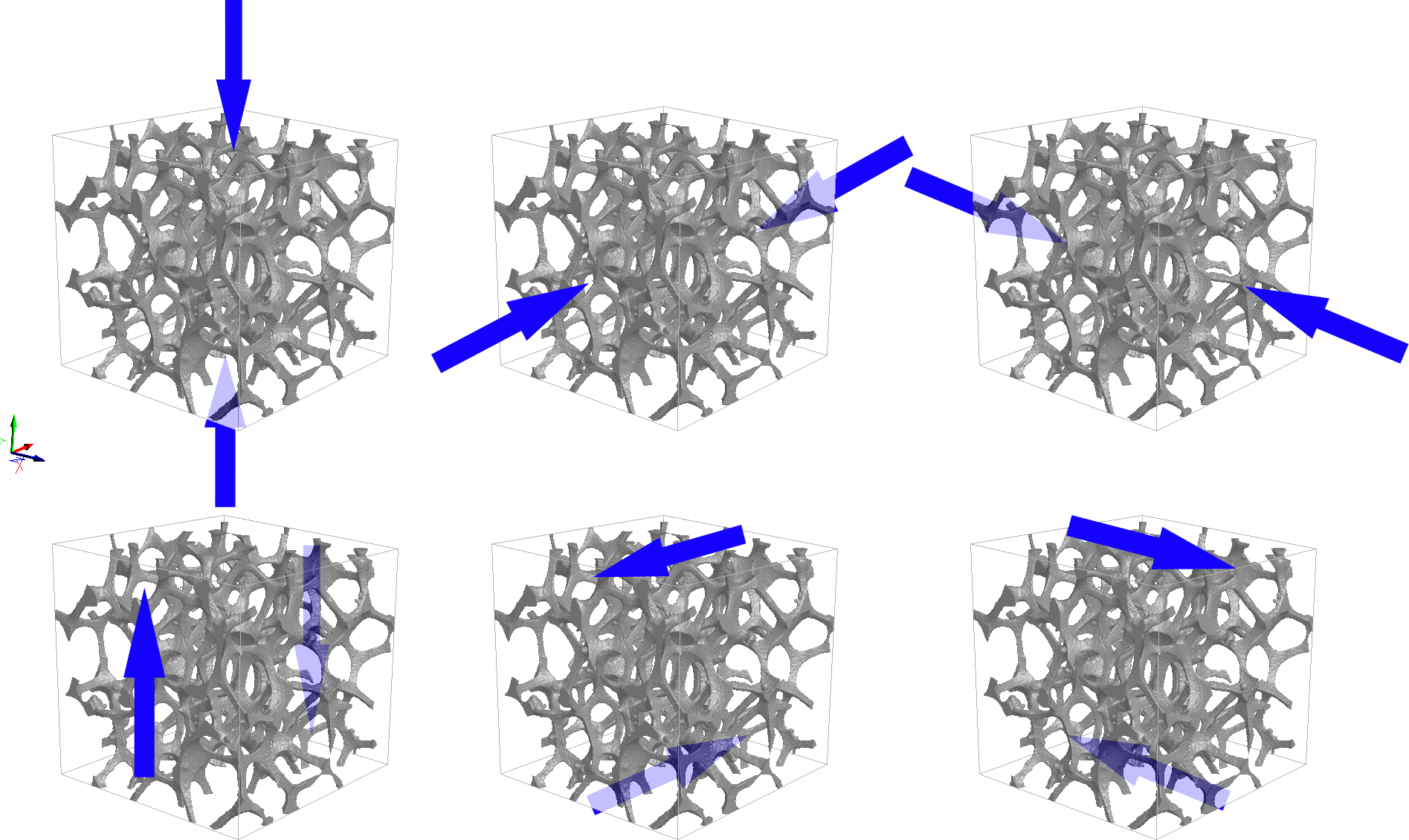}
\caption{Application of boundary conditions for homogenization. Blue arrows indicate the forces applied.}
\label{fig:homog}
\end{figure}

The micro-structure is captured in volume element $\Omega\subseteq\R^3$ with volume $V(\Omega)$
and boundary $\partial \Omega$.
The effective stiffness of this structure $\mathbb C^\ast$ connects the averaged stresses $\langle \sigma \rangle$ and strains $\langle \varepsilon \rangle$ via
\begin{equation}
    \label{eq:staub_1}
   \langle \sigma \rangle= \mathbb C^\ast : \langle \varepsilon \rangle.
\end{equation}
The $:$ in equation \eqref{eq:staub_1} denotes the double inner product between the fourth order tensor $\mathbb C^\ast$ and the second order tensor $\langle \varepsilon \rangle$ and corresponds to the mapping of the strains to the stresses (compare Hooke's law in one-dimensional elasticity $\sigma=E\, \varepsilon$).
The averaged stresses and strains are defined as volume averages of their microscopic counterparts as
\begin{equation}
\label{eq:staub_2}
   \langle \sigma \rangle:= \frac{1}{V(\Omega)}\int\limits_\Omega \sigma(x) \, \mathrm d x \quad \mbox{and } \quad  \langle \varepsilon \rangle := \frac{1}{V(\Omega)} \int\limits_\Omega \varepsilon(x) \, \mathrm d x.
\end{equation}
For the numerical computation of the effective foam properties, a displacement, which corresponds to a constant strain in a homogeneous reference material, is applied to the boundary of the volume element.
Due to the possible anisotropy of the foam, six loadcases  are computed -- three tension, three shear as sketched in Figure~\ref{fig:homog}.

After solving the six microscopic boundary value problems, we compute the complete effective stiffness tensor $\mathbb C^\ast$.

\subsection{Lippmann-Schwinger Fast Fourier Transform based solver}
\label{sse:LSFFT}

In the following, details on the solution of the boundary value problem are outlined. Our microscopic simulation is based on the solution of the Lippmann-Schwinger (LS) equations for 
elasticity, see \cite{Lippmann1950} and \cite{Zeller1973}.
The required Green's operator is explicitly known in the Fourier space. Thus, these equations are solved efficiently by application of a fast Fourier transform (FFT).
%Therefore in the following the solution scheme is referred to as LS-FFT.
This LS-FFT scheme is implemented in ITWM's micro-structure solver FeelMath, see \cite{Kabel2013} and \cite{Kabel2016}, also available as module ElastoDict \cite{ElastoDict} in the commercial software GeoDict \cite{GeoDict}.
The LS-FFT solver enables precise computation of microscopic stresses and strains directly in $\micro$CT images or other voxel based three-dimensional structures.  It is applicable to porous structures and therefore suitable for the open aluminum foam considered here.

In the solver, the equilibrium equation of the Cauchy stress $\sigma$ 
\begin{equation}
\label{eq:staub_3}
    \mathrm {div } \, \sigma(x) = 0, \quad x \in \Omega
\end{equation}
is considered in the micro-domain $\Omega$.
The kinematics for the strains $\varepsilon$ depending on the displacements $u$ and the fluctuations $v$ read 
\begin{equation}
\label{eq:staub_4}
   \left. \begin{array}{lll}
         \varepsilon(u)(x)&=&\langle \varepsilon \rangle+\varepsilon(v)(x)  \\
          \varepsilon(v)(x)&=&\frac{1}{2}\left( \mathrm{grad} \, v(x) +\mathrm{grad}^t \, v(x) \right)
    \end{array} \right\} \quad x \in \Omega
\end{equation}
in terms of the applied macroscopic strain $\varepsilon$.
At the boundary of the micro-structure $\partial \Omega$ (anti-)periodic boundary conditions are applied via
\begin{equation}
\label{eq:staub_5}
    \left. \begin{array}{lr}
        v(x) & \qquad \#  \\
        \sigma(x) \cdot n(x) & -\# 
    \end{array} \right\} \quad x \in \partial \Omega.
\end{equation}
The symbol $\#$ refers to periodicity, i.e. the fluctuations at opposite faces of the boundary are equal, whereas $-\#$ denotes antiperiodicity, i.e. the tractions $\sigma \cdot n$ at opposite faces point into opposite directions, but have the same magnitude.
The set of underlying equations is completed by a constitutive equation for the microscopic constituents.
For the aluminum we restrict ourselves to the linear elastic case.
Thus the microscopic stresses and strains are connected by an elasticity tensor $\mathbb C$ via
\begin{equation}
\label{eq:staub_6}
    \sigma(x)= \mathbb C(x) : \varepsilon(x)
\end{equation}
which only depends on Young's modulus $E$ and Poisson's ratio $\nu$ of the aluminum alloy. 
In the pores, the stiffness is set to zero.
In general the presented approach is however suitable to capture more complex material behavior like inelasticity and rate-dependency, see  \cite{Kabel2017,Schneider2021,Staub2018}.

Next, we focus on reformulating the periodic boundary value problem into an integral expression of the Lippmann-Schwinger type.
To this end, we introduce a constant homogeneous reference stiffness tensor $\mathbb C^0$ instead of the stiffnesses of the aluminum and the pores.
This homogeneous stiffness tensor is applied to define the polarization tensor $\tau$ as
\begin{equation}
\label{eq:staub_7}
    \tau(x)=\sigma(x) - \mathbb C^0 : \varepsilon(x).
\end{equation}
With the help of Green's operator $\Gamma^0$ associated with the reference stiffness, the solution of the equilibrium equation \eqref{eq:staub_1} reads
\begin{equation}
    \label{eq:staub_8}
    \varepsilon(x)=\langle \varepsilon \rangle- \left(\Gamma^0 \ast \tau \right)(x).
\end{equation}
The convolution operator $\ast$ is defined by
\begin{equation}
    \label{eq:staub_9}
    \left(\Gamma^0 \ast \tau \right)(x)= \int\limits_\Omega \Gamma^0(x-y):\tau(y) \; \mathrm d y.
\end{equation}
Combining the constitutive law \eqref{eq:staub_4}, the definition of the polarization stress \eqref{eq:staub_7}, and the solution \eqref{eq:staub_8} yields the LS equation as
\begin{equation}
    \label{eq:staub_10}
    \varepsilon = \varepsilon(x)+ \Gamma^0(x) \ast \left( \mathbb (C(x)-\mathbb C^0):\varepsilon(x) \right) = \left(I+B_\varepsilon(x) \right) \varepsilon(x).
    \end{equation}
Green's operator $\Gamma^0$ does not depend on the fluctuations and thus only depends on the homogeneous reference stiffness $\mathbb C^0$, see \cite{Kroener1977}.

The LS equation \eqref{eq:staub_10} can be solved iteratively using the Neumann series expansion or by using the conjugate gradient method.
As Green's operator is explicitely known in Fourier space, the Fourier transform is applied. 
A more detailed description of the algorithm can be found in \cite{Staub2018}.

Note that the discretization by Fourier polynomials as presented in \cite{Moulinec1998} leads to convergence problems for porous structures due to the infinite stiffness contrast of the microscopic constituents.
Therefore we use a finite difference discretization based on a staggered grid, which converges also for highly porous materials, see \cite{Schneider2016b}.
A comprehensive overview of FFT-based homogenization methods is given in \cite{Schneider2021}.

\section{Application example}
\label{sec:application}

Now, all methods described above are applied to a real world open aluminum foam sample. 

\subsection{Material}

We consider one of the open-cell aluminium alloy foam samples investigated 
in \cite{jung2020}. The sample made by CellTec
Materials GmbH, Dresden, Germany consists of AlSi$_7$Mg$_{0.3}$, has nominal pore 
size 10\,ppi (pores per inch) and mean density 0.156\,g/cm$^3$
corresponding to a porosity of 94.2\,\%. 

\begin{figure}[b]
\sidecaption
\includegraphics[width=.4\textwidth]{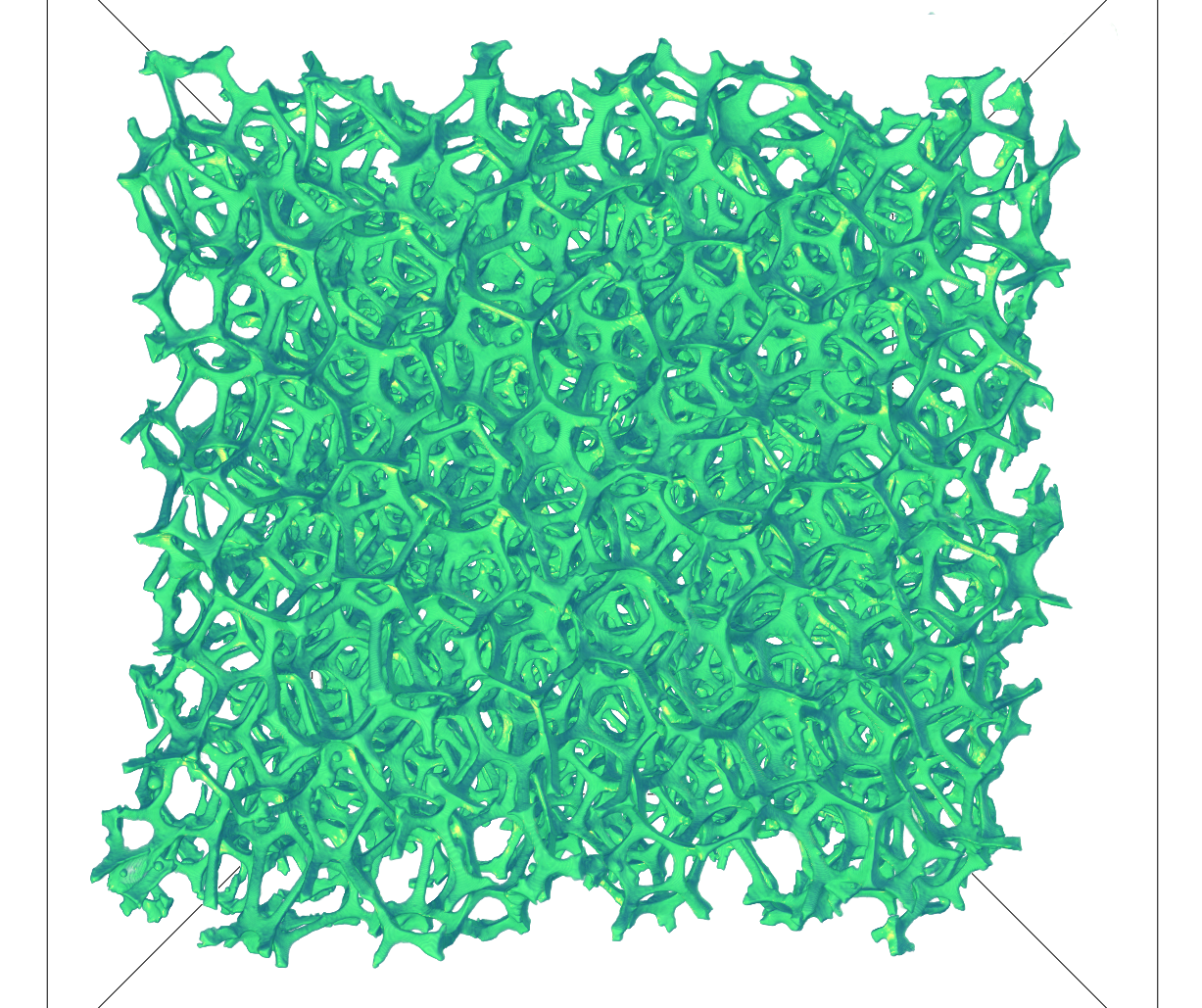}\
\includegraphics[width=.2\textwidth]{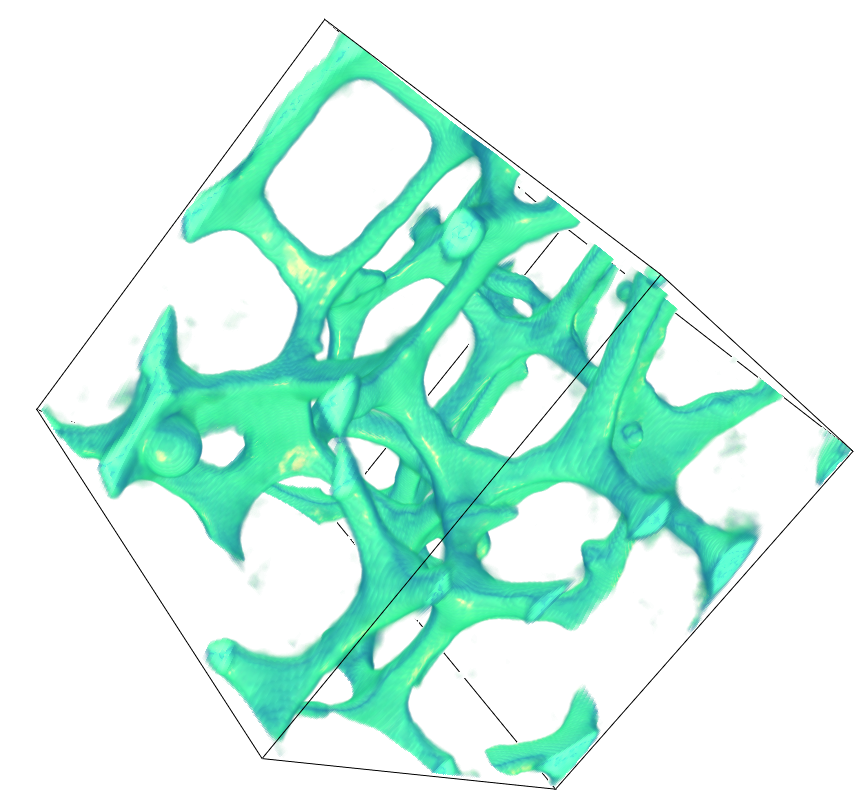}\
\includegraphics[width=.3\textwidth]{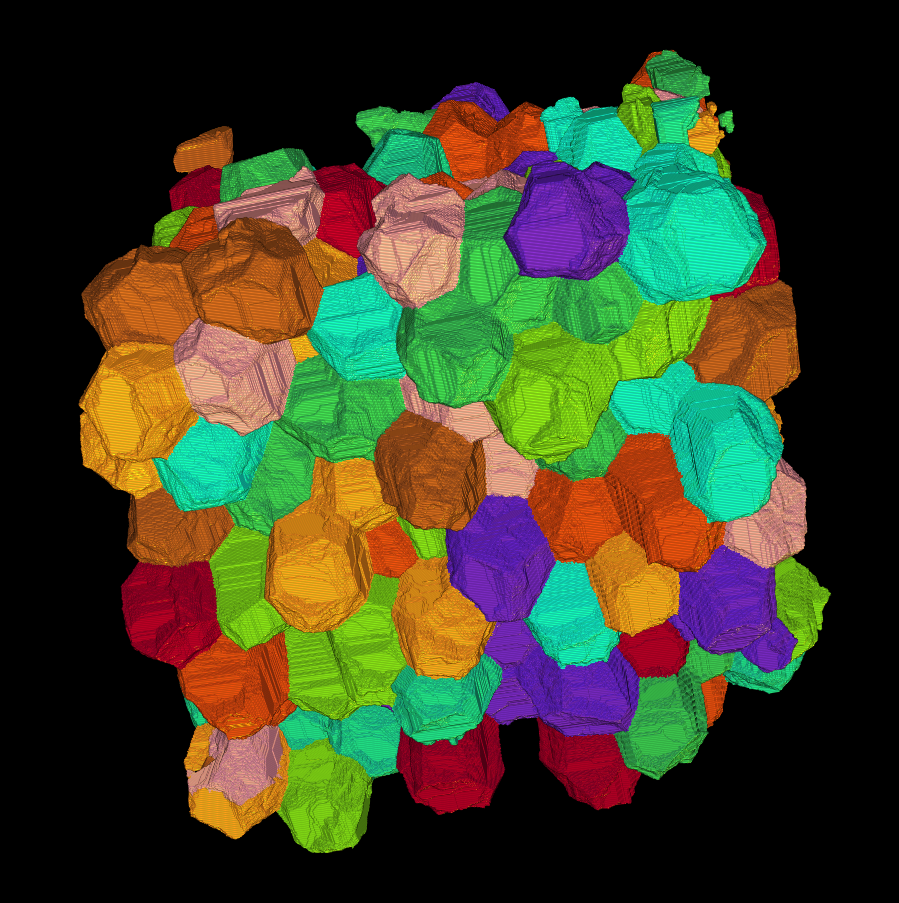}
\caption{Left: Volume rendering of the reconstructed CT image. The CT image taken at ITWM has originally a voxel size of $29.9\,\micro$m and the sample is contained in a cube of edge length 1\,500 voxels corresponding to 4.5\,cm. Center: sub-volume of edge length 0.9 cm with blob like production leftover. Right: System of reconstructed cells.}
\label{fig:vr-full} 
\end{figure}

A cubic sample of edge length 40\,mm is spatially imaged by $\micro$CT at voxel size $29.44\,\micro$m. See Figure~\ref{fig:vr-full} for a volume rendering.

Larger samples of size 40\,mm $\times$ 40\,mm $\times$ 80\,mm were tested mechanically. 
Uniaxial compression and tensile tests were performed, where the samples had to 
be infiltrated by a resin to allow for clamping for the tensile tests. 
Details are described in \cite{Jung2017}. In \cite{jung2020}, the foam 
structure was image analytically separated into vertices and struts. 
Mechanical behaviour of struts was 
investigated for five classes divided according to orientation.  

\subsection{Image analysis and model fit}
\label{ssec:modelFit}

The CT image of the foam was binarized by a global grayvalue threshold chosen according to Otsu's method. Subsequently, the cells were reconstructed by the watershed on iEDT approach described in Section ~\ref{ssec:segmentation}. Due to contamination of the foam structure (see Figure~\ref{fig:vr-full}, middle), some cells were erroneously split into two parts. These errors were corrected manually. 

Estimation of the mean cell diameters in the three coordinate directions reveals that cells are elongated by a factor $s=1.2526$ in $y$-direction. 
The model fitting procedure detailed in Section~\ref{subsec:modelfitting} is formulated for isotropic structures. 
Hence, the structure is scaled by $1/s$ prior to the fitting. 
The cell characteristics estimated from the isotropic structure are listed in Table~\ref{tab:Al_foam}. 
Cells intersecting the boundary of the image are not observed completely. 
Hence, they are not included in the statistics. 
As larger cells are more likely to intersect the image boundary, ignoring boundary cells will result in a sampling bias. 
This is corrected by weighting the cells according to the inverse of their probability of being observed (Miles-Lantuejoul correction, \cite[pp.~246]{ser82}). 
In total, 336 cells with a total volume of $26\,246.11 \text{mm}^3$ were included in the statistics. 

The model fitting procedure returns a packing fraction of $\kappa =60 \%$ and a coefficient of variation $c= 0.414$ as optimal parameters when assuming a lognormal distribution of sphere volumes. Using these parameters, realizations of the tessellation models containing 820 cells are simulated in the unit cube. After rescaling the $y$-axis by the factor $s$, the tessellation edges are discretized into volume images. The voxel spacing for the discretization is chosen identical to the CT image, i.e. $29.44 \, \micro$m. The images are cropped to a size of $1\,359^3$ voxels to obtain the original sample edge length of $40$ mm. 

\begin{figure}[b]
\sidecaption
\includegraphics[width=.24\textwidth]{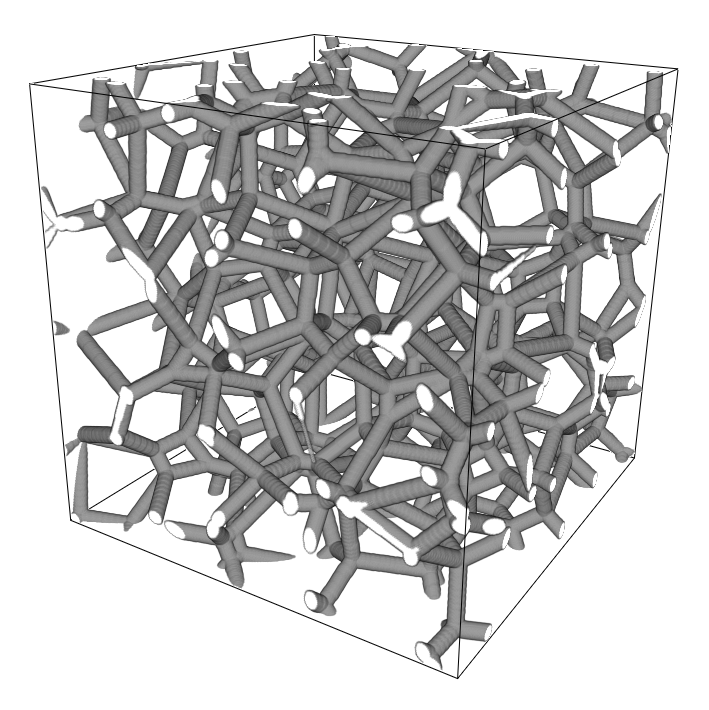}
\hfill
\includegraphics[width=.24\textwidth]{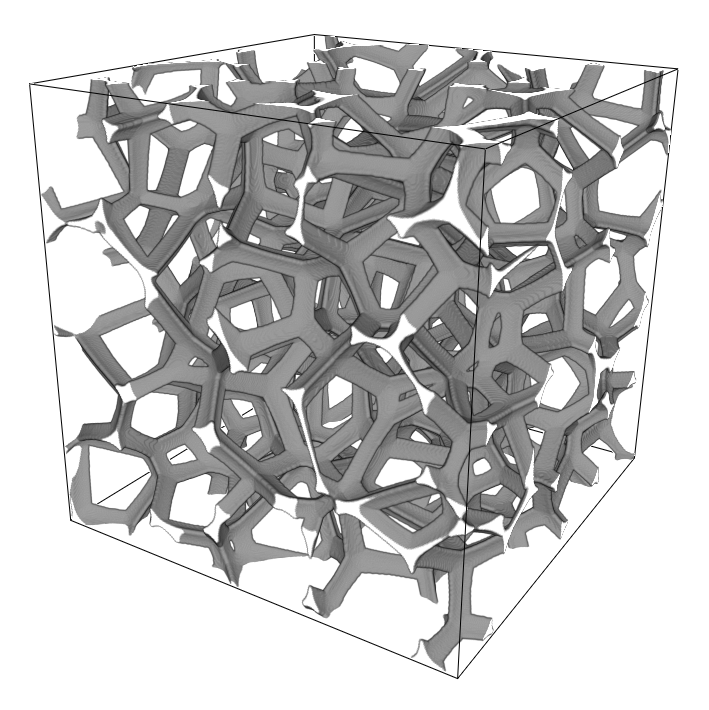}
\hfill
\includegraphics[width=.24\textwidth]{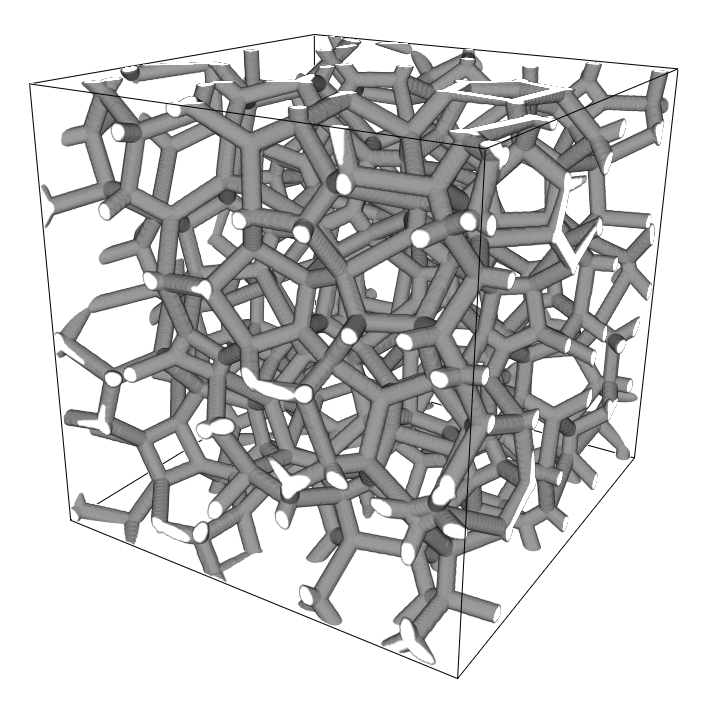}
\hfill
\includegraphics[width=.24\textwidth]{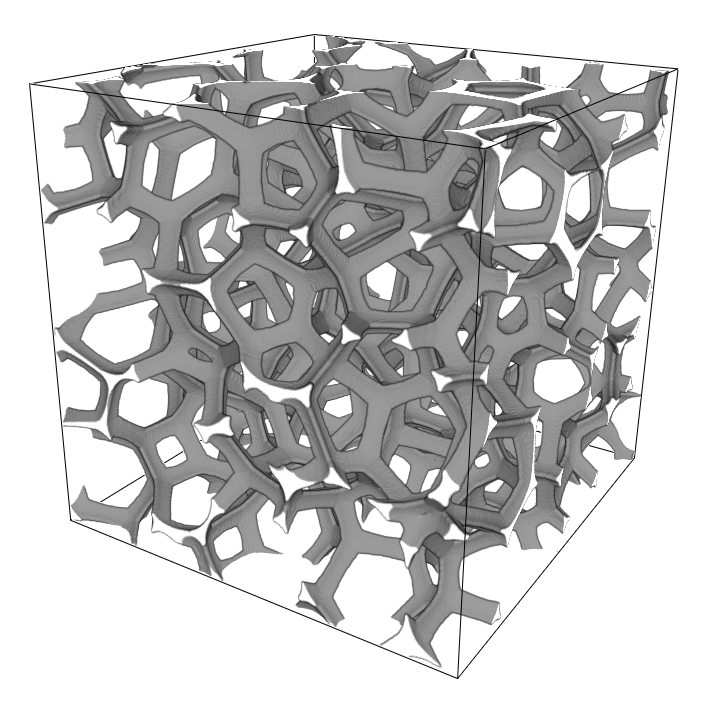}
\caption{Volume renderings of four synthetic foams derived from the same Laguerre tessellation model realization. From left to right: circular and triangular strut cross sections and their relaxed versions. Sub-volumes of $600^3$ voxels.}
\label{fig:vis-model}
\end{figure}

Dilation of the edge system of the tessellations finally yields synthetic foam structures. Strut thickness is chosen such that the volume fraction is close to the value $V_V=5.755\%$ estimated from the binarized CT image of the real foam. 
We consider two cross-sectional shapes of the struts -- the simplest case of 
perfectly circular strut cross sections of constant radius and concave triangular struts. 
Additionally, to investigate the effect of relaxation, 
the tessellations are relaxed by using the Surface Evolver 
and also discretized with the two choices of strut cross sectional shape. That way, 
each realization of the Laguerre tessellation model yields four synthetic foam 
structures. All four synthetic foams are visualized in Figure~\ref{fig:vis-model}.

\begin{table}[!htp]
\caption{Estimated mean values and standard deviations of the cell characteristics of the aluminium foam and the best fit models for lognormally distributed volumes of the generating spheres.}
\label{tab:Al_foam}
\begin{center}
\begin{tabular}[t]{l r r r r r r r  r} 
\hline\noalign{\smallskip}
& & \multicolumn{2}{c}{scaled data}& & \multicolumn{4}{c}{isotropic model}\\
& & mean &std &&mean& deviation&std& deviation\\
\noalign{\smallskip}\svhline\noalign{\smallskip}
$v$ & $[{\rm mm}^3]$& 79.738&21.972&$\quad$&79.739&$\pm$0.00\% & 26.598&$\pm$21.06\%\\
$s$ & $[{\rm mm}^2]$&104.468&22.241&$\quad$&98.608&$-5.61$\%&19.832&$-10.83$\%\\
$d$ & $[{\rm mm}]$&5.811&0.740&$\quad$&5.806&$-0.08$\%&0.568&$-23.22$\%\\
$F_C$& &13.828&2.256&$\quad$&14.06&$+1.69$\%&2.173&$-3.71$\%\\
\noalign{\smallskip}\hline\noalign{\smallskip}

\end{tabular}
\end{center}
\end{table}

%\clearpage

\subsection{Prediction of mechanical properties}
\label{subsec:mechprops}
We now predict the mechanical properties of the aluminum foam based on the binarized CT image and on realizations of the four synthetic foams.

\subsubsection{Validation of simulation}
\label{ssubsec:validation}
Before actually predicting, we validate the numerical computation of the effective quantities with ElastoDict \cite{ElastoDict} by comparing with experimental data.
To this end, a compression load case with a prescribed strain is simulated directly on the binarized CT image of the foam. %, for details see section \ref{ssec:modelFit}.
The resulting effective stresses and the corresponding stiffness in the loading direction are calculated and compared to those determined by a compression test on the foam as described in \cite{Jung2017}.

First, we summarize the solver settings: The considered foam structure is not periodic. Hence, symmetric boundary conditions with free deformations in tangential directions are chosen.
These symmetric boundary conditions and their compatibility with the periodic boundary conditions from Section~\ref{sse:LSFFT} are described in detail in \cite{Grimm-Strele2021}. Thanks to the compatibility, this type of boundary condition is often referred to as periodicity compatible mixed uniform boundary conditions (PCMUBC). Generating a periodic structure by means of mirroring and effective implementation of PCMUBC in the context of FFT are connected in \cite{Grimm-Strele2021}, too.
 
We apply composite voxels to reduce computation time as suggested in \cite{Kabel2015}. 
In this approach, the image is downsampled and the stiffness of the resulting mixed voxels or composite voxels is determined with the appropriate material of the rotated laminate, i.e. not only the local volume fraction is taken into account. 
Additionally, this mixing rule incorporates directional information on the material interface. 
Applying this strategy, we cut $1\,314\times1\,335\times1\,285$ voxels down to $328\times333\times321$ voxels, and consequently reduce CPU time and memory consumption significantly while keeping the loss of accuracy acceptable.
%Choosing the Memory Efficient CG scheme \cite{Grimm-Strele2019} further reduces the required memory.

Second, we care for the material parameters: We restrict ourselves to the linear elastic case here. 
Therefore only Young's modulus $E_{Al}$ and Poisson's ratio $\nu_{Al}$ of the aluminum at the micro-scale are required.
Poisson's ratio is taken from literature as $\nu_{Al}=0.33$.
Young's modulus of the aluminum is chosen based on the experimental data from \cite{Heinze2018} measured in fives different pores. The values are summarized in Table~\ref{tab:Young}.
\begin{table}[!htp]
\caption{Experimentally identified Young's moduli for the aluminum alloy \cite{Heinze2018}. }
\label{tab:Young}
\begin{center}
\begin{tabular}[t]{l c } 
\hline\noalign{\smallskip}
Pore no. & Young's modulus [MPa] \\
\noalign{\smallskip}\svhline\noalign{\smallskip}
1 & 4\,112 \\
2 & 3\,832 \\
3 & 3\,768 \\
4 & 4\,152 \\
5 & 3\,978 \\
\noalign{\smallskip}\hline\noalign{\smallskip}
\end{tabular}
\end{center}
\end{table} 
In the simulations, the average of these moduli (4\,000~MPa) and the mimimal and maximal values are considered. 
Attention should be paid to the fact that these values differ significantly from those usually reported for aluminum alloys ($\approx  70 \cdot 10^3$~MPa). These deviations are due to micro-porosity of the struts and inclusions resulting from the manufacturing process \cite{Luksch2021}.  

\begin{figure}[ht!]
\sidecaption
\includegraphics[width=0.99\textwidth]{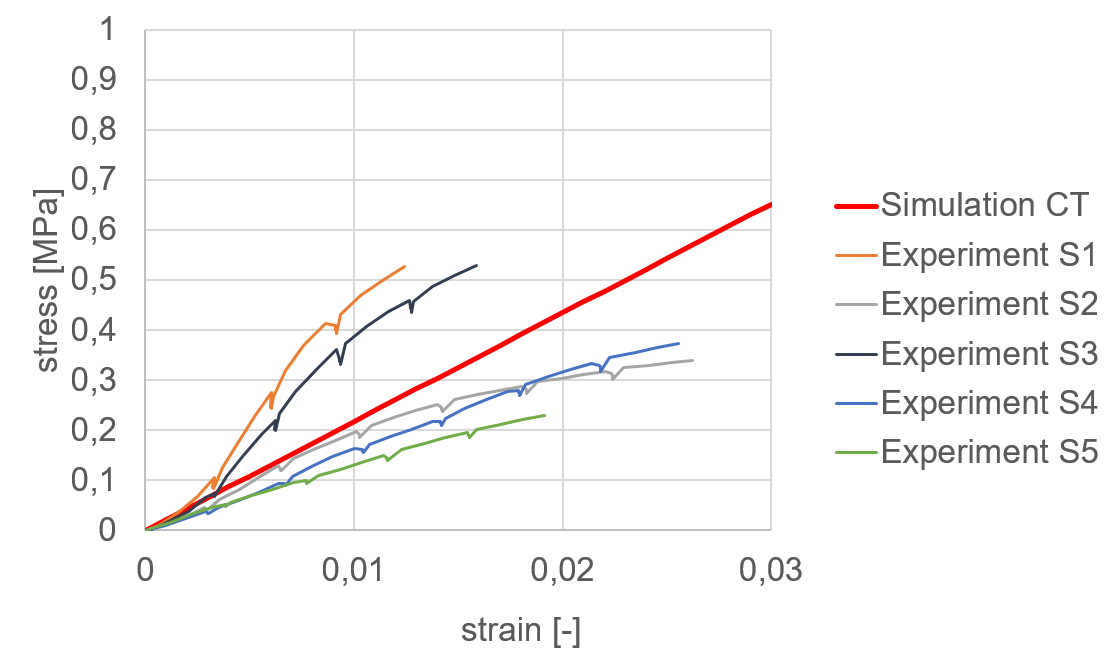}
\caption{Comparison of simulation and experiment in longitudinal direction $y$ in terms of stress-strain diagram, Young's modulus of the aluminum $E_{Al}=$4\,000~MPa.}
\label{fig:simExp1a}       
\end{figure}

Figure~\ref{fig:simExp1a} displays the experimental stress-strain curves of the five foam samples under compression in the linear loading regime, see \cite{Jung2017} for the complete loading and unloading curves,  and the effective stress-strain curves simulated directly in the complete binarized CT image. Simulated and experimental data match very well.
Young's modulus of the aluminum alloy is chosen as the average of the five samples summarized in Table~\ref{tab:Young}.

\begin{figure}[ht!]
\sidecaption
\includegraphics[width=0.6\textwidth]{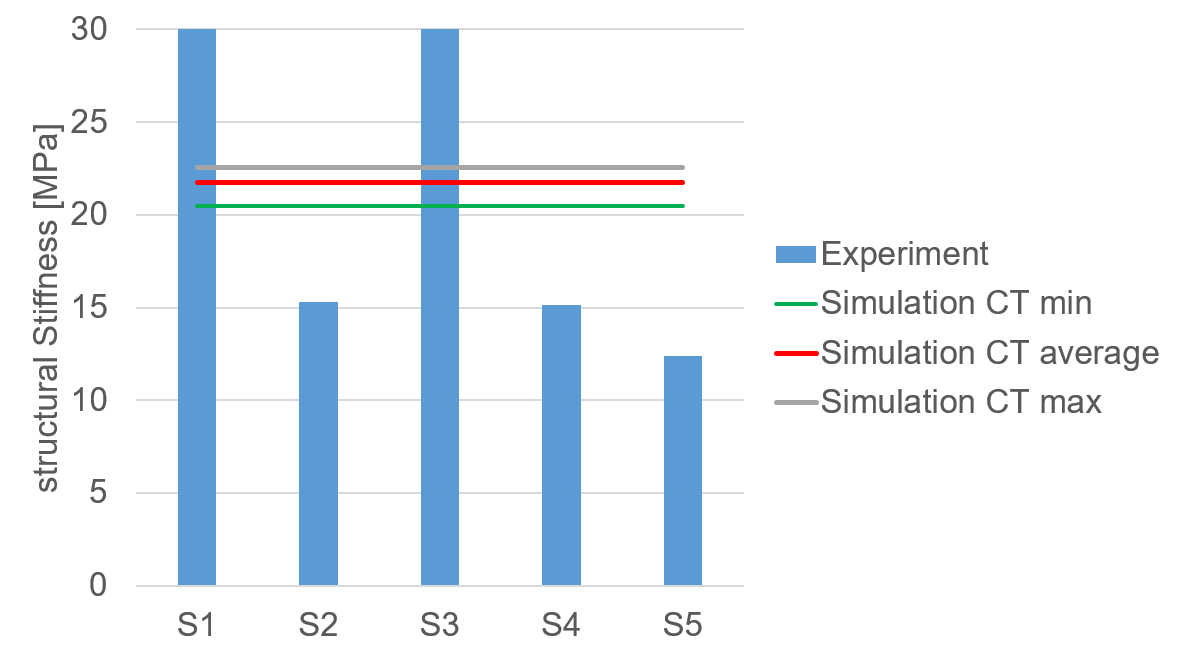}
\caption{Comparison of simulation and experiment in longitudinal direction $y$ in terms of structural stiffness for varying Young's moduli of the aluminum alloy.}
\label{fig:simExp1b}       
\end{figure}

Figure~\ref{fig:simExp1b} shows how the microscopic Young's modulus of the aluminum alloy influences the simulation result in terms of the effective modulus, i.e. the slope of the effective stress-strain curve. 
As input for the simulation, the minimal, maximal, and average values from Table~\ref{tab:Young} are considered. 
These simulation results are compared to the structural stiffness obtained by the compression experiments.

The structural stiffness corresponds to the slope of the experimental stress-strain curves in Figure~\ref{fig:simExp1a}. 
Clearly, variations of the microscopic Young's modulus in the considered range have only minor influence on the effective modulus.
Therefore, we take only the average value of 4\,000~MPa into account in the following.

\begin{figure}[ht!]
\sidecaption
\includegraphics[width=0.99\textwidth]{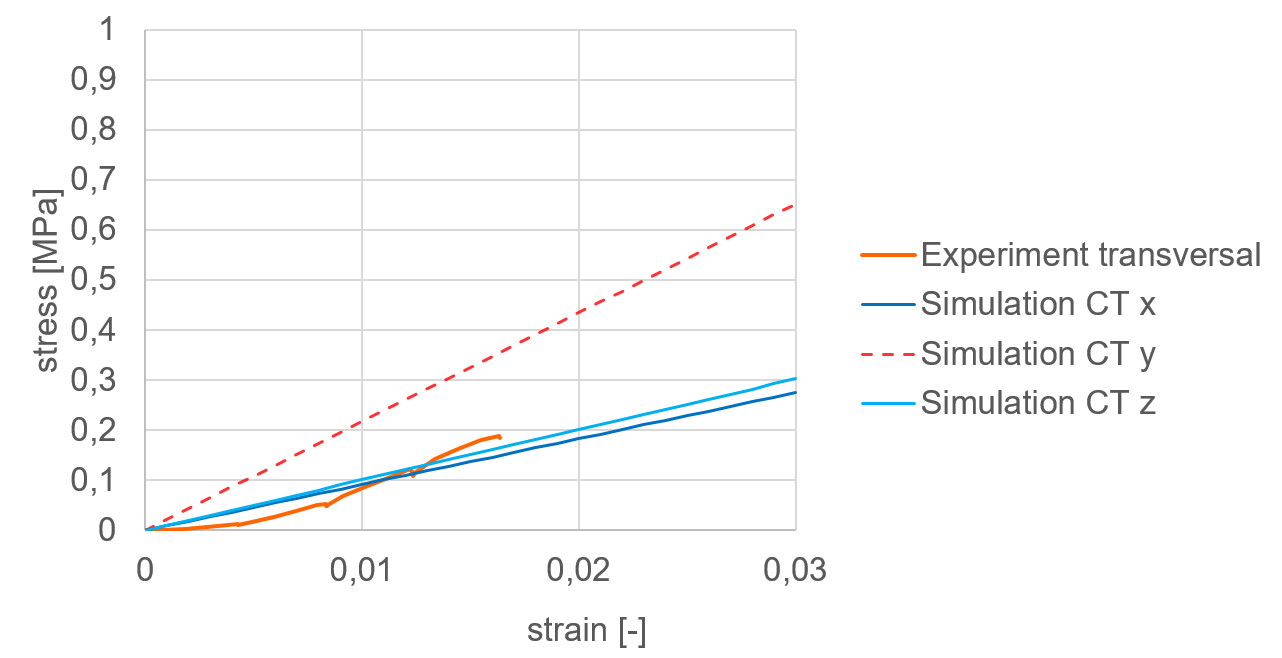}
\caption{Comparison of simulation and experiment in transversal direction, Young's modulus of the aluminum $E_{Al}=$4\,000~MPa.}
\label{fig:simExp2}       
\end{figure}

Validation of the simulation is completed by a comparison of the effective stresses simulated in the transversal directions (here $x$ and $z)$ and a compression experiment in the corresponding directions in Figure~\ref{fig:simExp2}. The simulated effective stresses in $x$- and $z$-directions are almost the same and they are much smaller than those in $y$-direction.
Thus, the foam is much stiffer in the longitudinal $y$-direction.
Figure~\ref{fig:simExp2} also shows that simulation and experiment also fit very well in the transversal directions. 

In the next section the validated simulation is applied to the synthetic foams derived from realizations of the Laguerre tessellation model fit in Section~\ref{ssec:modelFit}.

\subsubsection{Effective mechanical properties of the synthetic foams}
\label{ssubsec:prediction}

In order to predict the properties of the synthetic foams, five realizations of the Laguerre tessellation model yielding altogether 20 synthetic foams with circular or concave triangular struts, in relaxed state or not, are considered.
All load cases displayed in Figure~\ref{fig:homog} are applied to each synthetic foam and the full effective stiffness tensor is computed.
The effective Young's modulus in each direction ($E_x, E_y$ and $E_Z$) is approximated orthotropically, i.e. it is assumed that the inverse effective stiffness tensor reads
\begin{equation}
    \label{eq:ortho}
    \mathbb C^{-1} = \left(\begin{array}{cccccc}
         \frac{1}{E_x}& - \frac{\nu_{yx}}{E_y} & - \frac{\nu_{zx}}{E_z} & 0 & 0 &0  \\
          -\frac{\nu_{xy}}{E_x} &  \frac{1}{E_y} &  -\frac{\nu_{zy}}{E_z} & 0  &0 & 0 \\
          - \frac{\nu_{xz}}{E_x}  -  \frac{\nu_{yz}}{E_y} & \frac{1}{E_z} & 0 & 0 & 0 \\
          0 & 0 & 0 & \frac{1}{G_{yz}} & 0 & 0 \\
          0 & 0 & 0 & 0 & \frac{1}{G_{zx}} & 0 \\
          0 &  0 & 0 & 0 & 0 & \frac{1}{G_{xy}}
    \end{array} \right),
\end{equation}
where $\nu_{xy}, \nu_{xz},...$ denote Poisson's ratios and $G_{xy}, G_{zx}, G_{xy}$ the shear moduli in the corresponding directions.

\begin{figure}[ht!]
\sidecaption
\centering
\includegraphics[width=0.6\textwidth]{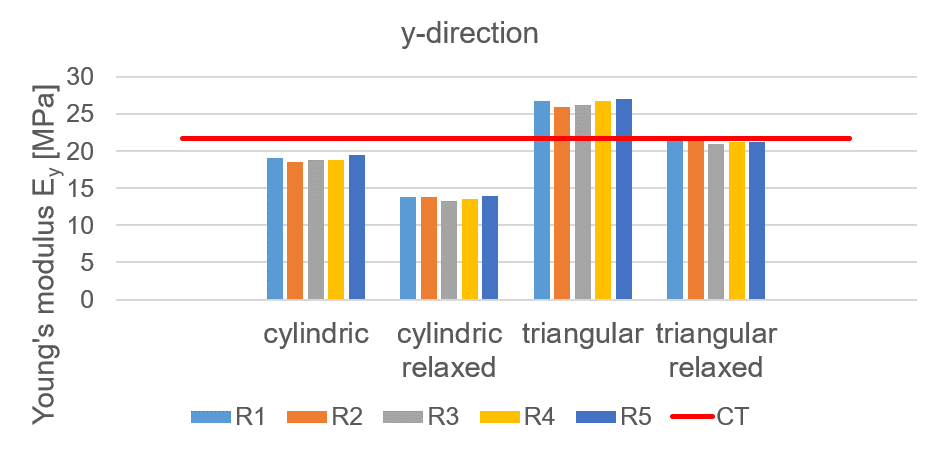}
\caption{Comparison of effective stiffness of synthetic and real foams in longitudinal direction. CT is the real foam's micro-structure as given by the binarized CT image. R1 to R5 are the five realizations of the fitted Laguerre tessellation model.}
\label{fig:modely}       
\end{figure}

We compare our simulation results to those obtained by direct simulations in the CT image. For the longitudinal ($y$) direction, the variations within one synthetic foam type do not influence the effective stiffness significantly, see Figure~\ref{fig:modely}.
For the two synthetic foam models with circular struts, the stiffness is underestimated.
On the other hand, concave triangular struts lead to overestimation of the stiffness. The true strut cross section observed in the CT image of the real foam is in between both models - it is clearly triangular, but less concave than the model.
The relaxed version of the synthetic foam with triangular struts is nevertheless able to reproduce the stiffness behavior in longitudinal direction of the real foam very well as the relaxation procedure reduces the effective stiffness for both strut geometries.

\begin{figure}[ht!]
\sidecaption
\includegraphics[width=0.5\textwidth]{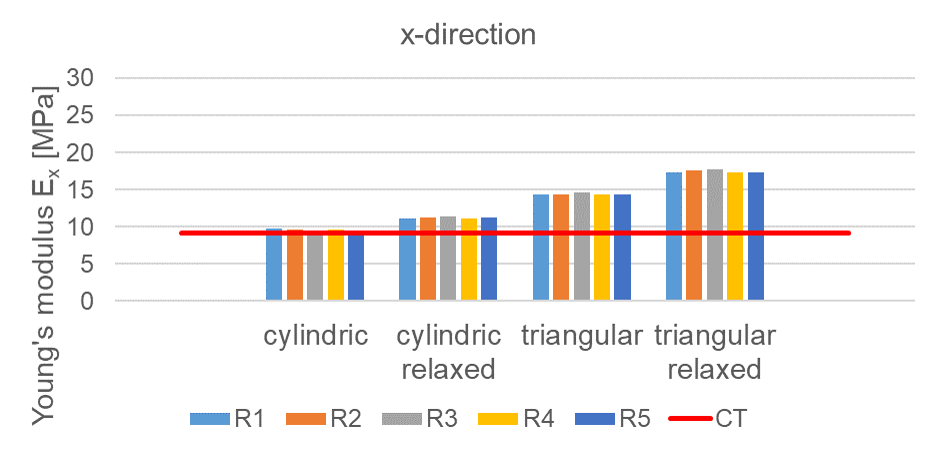}
\includegraphics[width=0.5\textwidth]{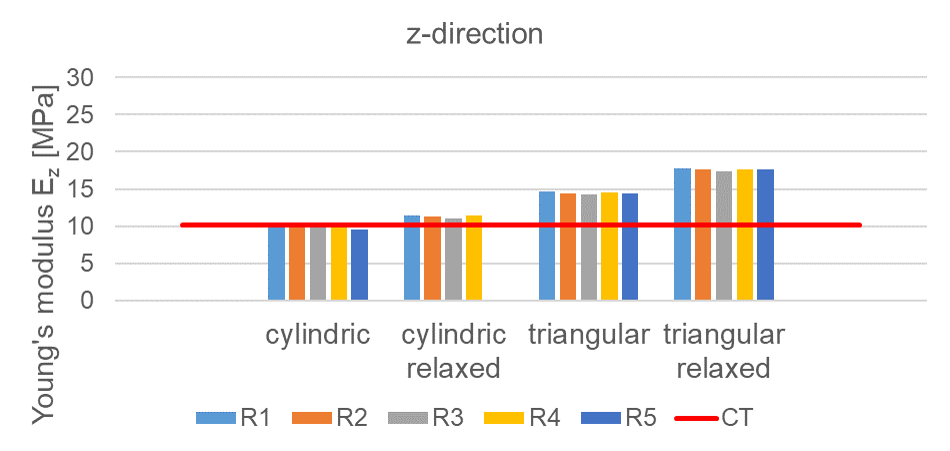}
\caption{Comparison of effective stiffness of synthetic (derived from R1 to R5) and real (CT) foams in transversal directions.}
\label{fig:modelxz}       
\end{figure}

Analysis of the mechanical properties of the synthetic foams is completed by the stiffness in the transversal directions ($x$ and $z$) reported in Figure~\ref{fig:modelxz}.
The anisotropic behavior of the real foam is captured by all synthetic foam models.
The stiffness in the transversal directions is smaller than in longitudinal direction, agreeing with the experimental results as well as with those obtained by direct simulation in the binarized CT image.
Not surprisingly, the two relaxed synthetic foams behave less anisotropic than their unrelaxed counterparts.
Relaxation increases the effective stiffness in transversal direction, while it decreases  the effective stiffness in  longitudinal direction.
The unrelaxed synthetic foam with the circular struts matches the mechanical behavior of the foam in transversal direction best.

In summary, all four synthetic foam models reproduce the mechanical compression 
behavior of the real aluminum foam well, including the anisotropy.
Taking into account all three directions, the effective stiffness of the real foam is best reproduced by the 
unrelaxed synthetic foam with circular strut cross-sections, although the 
stiffness in longitudinal direction is slightly underestimated.

\section{Conclusion}

In this contribution, we describe completely the well-established, nevertheless
rarely fully expanded workflow for simulating mechanical properties of a material 
based on stochastic geometry models of its micro-structure. 
Our use case is an open metal foam. 
Consequently, a random tessellation model is fit to it.

Stochastic geometry models not only capture naturally the microscopic heterogeneity 
of materials structures. Fitting them to the observed structure allows to 
generate many realizations as large as needed. Moreover, the effects of selected 
micro-structural geometric features can be investigated independently 
as demonstrated here for cross-sectional shape of the struts and relaxation. 
Model fitting and careful validation of the simulations further open the 
opportunity to determine the size of the representative volume and to 
devise optimized micro-structures.

Validation or calibration of the simulations by experiments is 
indispensable as the respective properties of the bulk material needed as 
input are hard to get. 
The values for the material e.~g. in the struts can differ significantly from 
tabular values due to effects of the production process on scales finer than the 
micro-structure. 

The elastic properties of the aluminum foam could be reproduced well by all  
four synthetic foams considered here. 
More sophisticated geometries might be needed when moving on to plastic properties. 
For instance, the local thickness of the struts could be fit
\cite{jung2020,liebscher2013}.  
Moreover, real 
cellular structures can feature much more complex anisotropies due to 
production processes being based on polymer foams. In the latter, the foaming direction usually stands out. The structure might however be additionally distorted e.~g. by the movement of a conveyor belt during foaming. Adequate characterization methods and modelling are subject of current research.

\begin{acknowledgement}
We thank Franz Schreiber for the CT scans, Felix Arnold for assisting in the image processing of the CT image, and Sonja F\"ohst for providing the code for model simulation.
\end{acknowledgement}

%\bibliographystyle{spmpsci}
%\bibliography{litbank}

\end{document}